\documentclass[11pt]{article}

\usepackage{amsfonts,amsmath}
\newtheorem{theorem}{Theorem}[section]

\newtheorem{example}{Example}
\newtheorem{assumption}{Assumption}
\newtheorem{cor}{Corollary}[section]

\numberwithin{equation}{section}
\def\la{\langle}
\def\ra{\rangle}

\def\x{\mathbf{x}}

\def\R{\mathbb{R}}

\def\N{\mathbb{N}}
\def\K{\mathbf{K}}

\def\Q{\mathbf{Q}}

\def\P{\mathbf{P}}

\def\A{\mathbf{A}}
\def\V{\mathbf{V}}

\def\B{\mathbf{B}}
\def\c{\mathbf{C}}

\def\f{\mathbf{f}}

\def\v{\mathbf{v}}

\def\b{\mathcal{B}}
\def\a{\mathbf{a}}

\def\s{\mathcal{S}}

\def\c{\mathbf{c}}
\def\b{\mathbf{b}}
\def\s{\mathcal{S}}

\def\la{\langle}
\def\ra{\rangle}
\def\blambda{{\boldmath{\lambda}}}

\title{A Lagrangian relaxation view of linear and semidefinite hierarchies}

\author{Jean B. Lasserre\thanks{LAAS-CNRS and Institute of Mathematics,
       University of Toulouse,  7 Avenue du Colonel Roche,
       31077 Toulouse, France. {\tt lasserre@laas.fr}}}
\begin{document}
\date{}
\maketitle
\begin{abstract}
We consider the general polynomial optimization problem
$\P:\:f^*=\min \{f(\x)\,:\,\x\in\K\}$
where $\K$ is a compact basic semi-algebraic set.  We first show that the standard Lagrangian relaxation
yields a lower bound as close as desired to 
the global optimum $f^*$, {\em provided} that it is applied 
to a problem $\tilde{\P}$ equivalent to $\P$, in which sufficiently many
redundant constraints (products  of the initial ones) are added to the initial description of $\P$.
Next we show that the standard hierarchy of LP-relaxations of $\P$ (in the spirit of Sherali-Adams' RLT)
can be interpreted as a brute force simplification of the above Lagrangian relaxation
in which a nonnegative polynomial (with coefficients to be determined) is replaced with a constant polynomial equal to zero.
Inspired by this interpretation, we provide a systematic improvement of the LP-hierarchy
by doing a much less brutal simplification which results into a parametrized hierarchy of semidefinite programs
(and not linear programs any more).
For each semidefinite program in the parametrized hierarchy, the semidefinite constraint 
has a fixed size $O(n^k)$, independently of the rank in the hierarchy, in contrast with the standard hierarchy of semidefinite relaxations.
The parameter $k$
is to be decided by the user. When applied to a non trivial class of convex problems,
the first relaxation of the parametrized hierarchy is exact, in contrast with the 
LP-hierarchy where convergence cannot be finite. When applied to 0/1 programs
it is at least as good as the first one in the hierarchy of semidefinite relaxations.
However obstructions to exactness still exist and are briefly analyzed.
Finally, the standard semidefinite hierarchy can also be viewed as a simplification of
an extended Lagrangian relaxation, but different in spirit as sums of squares (and not scalars) multipliers are allowed.
\end{abstract}

\section{Introduction}

Recent years have seen the development of 
(global) semi-algebraic optimization and in particular 
LP- or semidefinite relaxations for the polynomial optimization problem:
\begin{equation}
\label{def-pb}
\P:\quad f^*=\displaystyle\min_\x\:\{f(\x)\::\: \x\in\K\:\}
\end{equation}
where $f\in\R[\x]$ is a polynomial and $\K\subset\R^n$ is the basic semi-algebraic set
\begin{equation}
\label{setk}
\K\,=\,\{\,\x\in\R^n\::\: g_j(\x)\,\geq\,0,\quad j=1,\ldots,m\},
\end{equation}
for some polynomials $g_j\in\R[\x]$, $j=1,\ldots,m$.\\

In particular, associated with $\P$ are two hierarchies of convex relaxations: \\

- {\it Semidefinite} relaxations based on Putinar's certificate of positivity on $\K$ \cite{Putinar93},
where the $d$-th convex relaxation of the hierarchy is a semidefinite program which solves
the optimization problem
\begin{equation}
\label{put-cert}
\gamma_d=\max_{t,\sigma_j}\:\{\,t\::\: f-t=\sigma_0+\sum_{j=1}^n \sigma_j\,g_j\}.
\end{equation}
The unknowns $\sigma_j$ are sums of squares polynomials with the degree bound constraint ${\rm degree}\,\sigma_jg_j\leq 2d$, $j=0,\ldots,m$,
and the expression  in (\ref{put-cert}) is a certificate of positivity on $\K$ for the polynomial $\x\mapsto f(\x)-t$.\\

- {\it LP}-relaxations based on Krivine-Stengle's certificate of positivity on $\K$ \cite{Krivine64a,Stengle74},
where the $d$-th convex relaxation of the hierarchy is a linear program which solves
the optimization problem
\begin{eqnarray}
\nonumber
\theta_d\,=\,\max_{\lambda\geq0,t}\:\left\{t\::\:
f-t\right.&=&\sum_{(\alpha,\beta)\in\N^{2m}_{d}}\lambda_{\alpha\beta}\,\left(\prod_{j=1}^mg_j^{\alpha_j}\right)\times\\
\label{kriv-cert}
&&\left.\left(\prod_{j=1}^m(1-g_j)^{\beta_j}\right)\right\},
\end{eqnarray}
where $\N^{2m}_{d}=\{(\alpha,\beta)\in\N^{2m}:\sum_j\alpha_j+\beta_j\leq d\}$. The unknown are $t$ and the nonnegative scalars $\lambda=(\lambda_{\alpha\beta})$, and
it is assumed that $0\leq g_j\leq 1$ on $\K$ (possibly after scaling) and the family $\{g_i,1-g_i\}$ generates the algebra $\R[\x]$ of polynomials.
Problem (\ref{kriv-cert}) is an LP because stating that the two polynomials in both sides of ``$=$" are equal yields linear constraints 
on the $\lambda_{\alpha\beta}$'s.
For instance, the LP-hierarchy from Sherali-Adams' RLT  \cite{SheraliAdams90} and their variants \cite{SheraliAdams99} are of this form.
See more details in \S \ref{sherali}.

In both cases, $(\gamma_d)$ and $(\theta_d)$, $d\in\N$, provide two monotone nondecreasing sequences of lower bounds 
on $f^*$ and if $\K$ is compact then both converge to $f^*$ as one let $d$ increase.
For more details as well as a comparison of such relaxations 
the interested reader is referred to e.g. 
Lasserre \cite{lasserre-imperial,Lasserre02a} and Laurent \cite{Laurent05}, as well as Chlamtac and Tulsiani \cite{tulsiani} 
for the impact of LP- and SDP-hierarchies on approximation algorithms in combinatorial optimization.

Of course, in principle,  one would much prefer to solve LP-relaxations rather than semidefinite relaxations (i.e. compute
$\theta_d$ rather than $\gamma_d$) because present
LP-software packages can solve problems with millions of variables and constraints, which is 
far from being the case for semidefinite solvers. And so the hierarchy (\ref{put-cert}) applies to problems of modest size only 
unless some sparsity or symmetry is taken into account in which case specialized variants can handle problems of much larger size.
However, on the other hand, the LP-relaxations (\ref{kriv-cert}) suffer from several serious 
theoretical and practical drawbacks. For instance, it has been shown in \cite{Lasserre02a,lasserre-imperial} that the LP-relaxations 
{\it cannot} be exact for most convex problems, i.e.,  the sequence of the associated optimal values converges
to the global optimum only {\it asymptotically} and not in finitely many steps. Moreover, the LPs of the hierarchy are numerically ill-conditioned.
This is in contrast with the semidefinite relaxations (\ref{put-cert})
for which finite convergence takes place for convex problems where $\nabla^2f(\x^*)$ is positive definite at every minimizer $\x^*\in\K$ 
(see de Klerk and Laurent \cite[Corollary 3.3]{deklerk-convex}) and occurs at the first relaxation for SOS-convex\footnote{An SOS-convex polynomial is a convex polynomial whose Hessian factors as $L(\x)L(\x)^T$ for some rectangular matrix polynomial $L$.
For instance, separable convex polynomials are SOS-convex.} problems \cite[Theorem 3.3]{lasserre-convex}.
In fact, as demonstrated in recent  works of Marshall \cite{Marshall2} and Nie \cite{nie2012},
finite convergence is generic (even for non convex problems).

So {\em would it be possible to define a hierarchy of convex relaxations in between (\ref{put-cert}) and (\ref{kriv-cert}), i.e.,
with some of the nice features of the semidefinite relaxations but with a much less demanding computational effort
(hence closer to the LP-relaxations)}? 
This paper is a contribution in this direction.\\

{\bf Contribution.} This paper consists of two contributions:
In the first contribution which is of theoretical nature, we describe a new hierarchy of convex relaxations for $\P$ with the following 
feature. Each relaxation in the hierarchy is a finite-dimensional convex optimization problem of the form:
\begin{equation}
\label{newhierarchy-def}
\rho_d\,=\,\max_{\blambda}\:\{\:G_d(\lambda)\::\: \blambda\geq0\},
\end{equation}
where $G_d(\cdot)$ is the concave function defined by:

\begin{eqnarray}
\nonumber
G_d(\lambda)\,:=\,\min_\x\:\left\{f(\x)\right.&-&\sum_{(\alpha,\beta)\in\N^{2m}_{d}}\lambda_{\alpha\beta}\,\left(\prod_{j=1}^mg_j^{\alpha_j}(\x)\right)\times\\
&&\left.
\label{def-gedelambda}
\left(\prod_{j=1}^m(1-g_j(\x))^{\beta_j}\right)\:\right\}.\end{eqnarray}
Therefore $\rho_d\leq f^*$ for all $d$. And we prove that:

(a)  $\rho_d\geq\theta_d$ for all $d$, and so $\rho_d\to f^*$ as one let $d$ increase.

(b) For convex problems $\P$, i.e., when $f,-g_j$ are convex, $j=1,\ldots,m$, and Slater's condition holds, the convergence is finite and occurs at the first relaxation,
i.e., $\rho_1=f^*$, in contrast with the LP-relaxations (\ref{kriv-cert}) where convergence cannot be finite (and is very slow on simple trivial examples). In fact computing $\rho_1$ is just applying the standard dual method 
of multipliers (or Lagrangian relaxation) to the convex problem $\P$.

(c) For $0/1$ optimization, i.e., when $\K\subseteq\{0,1\}^n$, finite convergence takes place and the optimal value $\rho_d$ provides a better lower bound than 
the one obtained with Sherali-Adams' RLT hierarchy \cite{SheraliAdams90}. In fact, the latter is solving (\ref{kriv-cert}) with only a subset of the 
products that appear in (\ref{kriv-cert}).

(d) Finally, (\ref{newhierarchy-def}) has a nice interpretation in terms of the dual method of Non Linear Programming
(or Lagrangian relaxation). To see this, consider the optimization problem $\tilde{\P}_d$ defined by:
\begin{equation}
\label{pbtildeP}
\min_\x\:\{f(\x)\::\:g_j(\x)^{\alpha_j}(1-g_j(\x))^{\beta_j}\,\geq\,0,\:
(\alpha,\beta)\in\N^{2m}_{d}\}\end{equation}
which has same value $f^*$ as $\P$ because $\tilde{\P}_d$ is just $\P$ with additional redundant constraints;
and notice that $\tilde{\P}_1=\P$.
Then solving (\ref{newhierarchy-def}) is just applying the {\it dual method} of multipliers in Non Linear Programming to $\tilde{\P}_d$; see e.g. \cite[Chapter 8]{bertsekas}. In general one obtains only a lower bound
on the optimal value of $\tilde{\P}_d$ when $\P$ is not a convex program).
And so our result states that the Lagrangian relaxation applied to $\tilde{\P}_d$ provides a lower bound as close to $f^*$ as desired, provided that
$d$ is sufficiently large, i.e., {\em provided} that sufficiently many redundant constraints are added to the description of $\P$.

Note in passing that this provides a rigorous {\em rationale} for the well-known fact that
adding redundant constraints helps for solving $\P$. Indeed, even though the new problems $\tilde{\P}_d$, $d\in\N$,
are all equivalent to $\P$, their Lagrangian relaxations are {\it not} equivalent to that of $\P$.

\subsection*{Practical and computational considerations}

Our second contribution has a practical and algorithmic flavor. Even though (\ref{newhierarchy-def}) is a convex optimization problem,
evaluating $G_d(\lambda)$ at a point $\lambda\geq0$ requires computing the unconstrained global minimum
of the function
\begin{eqnarray}
\nonumber
\x\mapsto L_d(\x,\lambda)&:=& f(\x)-\sum_{(\alpha,\beta)\in\N^{2m}_{d}}\lambda_{\alpha\beta}\,\left(\prod_{j=1}^mg_j(\x)^{\alpha_j}\right)
\times\\
\label{functionL}
&&\left(\prod_{j=1}^m(1-g_j(\x))^{\beta_j}\right),\end{eqnarray}
an NP-hard problem in general. After all, in principle the goal of Lagrangian relaxation is to end up with a problem which is easier to solve than $\P$,
and so, in this respect, the hierarchy (\ref{newhierarchy-def}) is {\it not} practical. 

So in this second part of the paper, we first show that the LP-relaxations (\ref{kriv-cert})
can be interpreted as a way to ``restrict"  and simplify the hierarchy (\ref{newhierarchy-def}) by a simple and brute force trick, so as to make it tractable
(but of course less efficient). Namely, a certain nonnegative polynomial (whose coefficients have to be determined) is 
imposed to be the constant polynomial equal to zero! More precisely, the nonnegative vector $\lambda$ in (\ref{newhierarchy-def}) is
restricted to a polytope so as to make the polynomial $L_d$ in (\ref{functionL}) constant!
In fact, if one had initially defined the LP-relaxations (\ref{kriv-cert})
as this brute force (and even brutal) simplification of (\ref{newhierarchy-def}), it would have been hard to justify.

Inspired by this interpretation, we propose a systematic way
to define improved versions of the LP-hierarchy (\ref{kriv-cert})  by simplifying (\ref{newhierarchy-def}) in a much less brutal manner. 
 We now impose the same
 nonnegative polynomial $L_d-t$ to be an SOS polynomial of fixed degree $2k$ (rather than the zero polynomial in (\ref{kriv-cert})).
The increase of complexity is completely controlled by the parameter $k\in\N$ and is chosen by the user.
That is, in the new resulting hierarchy (parametrized by $k$), each LP of the hierarchy (\ref{kriv-cert}) now becomes a semidefinite program but whose
size of the semidefiniteness constraint is fixed and equal to ${n+k\choose n}$, independently
of the rank $d$ in the hierarchy.
(It is known that crucial for solving semidefinite programs is 
the size of the LMIs involved rather than the number of variables.)
The level $k=0$ of complexity corresponds
to the original LP-relaxations (\ref{kriv-cert}), the level $k=1$ corresponds to a hierarchy of semidefinite programs 
with an Linear Matrix Inequality (LMI) of size $(n+1)$, etc.  To fix ideas, 
let us mention that for $k=1$, the first relaxation  (i.e., $d=1$) is even stronger than the first relaxation of the hierarchy (\ref{put-cert})
as it takes into account products of linear constraints; and so for instance,
when applied to the celebrated MAXCUT problem, the first relaxation has the Goemans-Williamson's performance guarantee.
Moereover, when $k=1$ one obtains the so-called ``Sherali-Adams + SDP" hierarchy already used for
approximating some 0/1 optimization problems.

So an important issue is: {\it What do we gain by this increase of complexity?}\\

Of course, from a computational complexity point of view, one way got evaluate the efficiency of 
those relaxations is to analyze whether they help reduce integrality gaps, e.g. for some
0/1 optimization problems. For the level $k=1$ (i.e. the ``Adams-Sherali + SDP hierarchy") some negative results in this direction have been provided
in Benabbas and Magen \cite{benabbas2}, and in  Benabbas et al. \cite{benabbas1}.

But in a  different point of view, we claim that a highly desirable property for a general purpose method (e.g., the hierarchies (\ref{put-cert}) or (\ref{kriv-cert}))
aiming at solving NP-hard optimization problems, is to behave ``efficiently" when 
applied to a class of problems considered relatively ``easy" to solve. Otherwise one might raise reasonable doubts on  its
efficiency for more difficult problems, not only in a worst-case sense but also  in ``average". Convex problems $\P$ as in (\ref{def-pb})-(\ref{setk}), i.e., when $f,-g_j$ are convex,
form the most natural class of problems which are considered easy to solve by some standard methods of Non Linear Programming;
see e.g. Ben-tal and Nemirovski \cite{Bental01}.
We have already proved that the hierarchy 
(\ref{put-cert}) somehow recognizes convexity. For instance, finite convergence takes places as soon as $\nabla^2f(\x^*)$
is positive definite at every global minimizer $\x^*\in\K$ (see deKlerk and Laurent \cite{deklerk-convex}); moreover, SOS-convex programs are solved at the first step of the hierarchy as shown in Lasserre \cite{lasserre-convex}.
On the other hand, the LP-hierarchy (\ref{kriv-cert}) behaves poorly on such problems as the convergence cannot be finite; see e.g. Lasserre \cite{lasserre-imperial,Lasserre02a}.

We prove that the gain by this (controlled) increase of complexity is precisely to permit finite convergence (and at the first step of the hierarchy)
for a non trivial class of convex problems. For instance with $k=1$ the resulting hierarchy of semidefinite programs
solves convex quadratic programs exactly at the first step of the hierarchy. And more generally,
for $k>1$, the first relaxation is exact for SOS-convex\footnote{A SOS-convex polynomial is such that its Hessian matrix is SOS, i.e., 
factors as $L(\x)L(\x)^T$ for some rectangular matrix polynomial $L$.}
 problems of degree at most $k$.  On the other hand, 
 we show that for non convex problems,  exactness at some relaxation in the hierarchy still implies restrictive
 conditions.
 
\section{Main result}

\subsection{Notation and definitions}
Let $\R[\x]$ be the ring of polynomials in the variables
$\x=(x_1,\ldots,x_n)$.
Denote by $\R[\x]_d\subset\R[\x]$ the vector space of
polynomials of degree at most $d$, which forms a vector space of dimension $s(d)={n+d\choose d}$, with e.g.,
the usual canonical basis $(\x^\alpha)$ of monomials.
Also, denote by $\Sigma[\x]\subset\R[\x]$ (resp. $\Sigma[\x]_d\subset\R[\x]_{2d}$)
the space of sums of squares (s.o.s.) polynomials (resp. s.o.s. polynomials of degree at most $2d$). 
If $f\in\R[\x]_d$, write
$f(\x)=\sum_{\alpha\in\N^n_d}f_\alpha \x^\alpha$ in the canonical basis and
denote by $\f=(f_\alpha)\in\R^{s(d)}$ its vector of coefficients. Finally, let $\s^n$ denote the space of 
$n\times n$ real symmetric matrices, with inner product $\la \A,\B\ra ={\rm trace}\,\A\B$, and where the notation
$\A\succeq0$ (resp. $\A\succ0$) stands for $\A$ is positive semidefinite.  With $g_0:=1$, the quadratic module $Q(g_1,\ldots,g_m)\subset\R[\x]$ generated by
polynomials $g_1,\ldots,g_m$, is defined by
\[Q(g_1,\ldots,g_m)\,:=\,\{\sum_{j=0}^m\sigma_j\,g_j\::\:\sigma_j\in\Sigma[\x]\}.\]
We briefly recall two important theorems by Putinar \cite{Putinar93} and Krivine-Stengle \cite{Krivine64a,Stengle74} respectively,
on the representation of polynomials positive on $\K$,
\begin{theorem}
Let $g_0=1$ and  $\K$ in (\ref{setk}) be compact.

(a) If the quadratic polynomial $\x\mapsto M-\Vert\x\Vert^2$
belongs to $Q(g_1,\ldots,g_m)$ and if $f\in\R[\x]$ is strictly positive on $\K$ then $f\in Q(g_1,\ldots,g_m)$.

(b) Assume that $0\leq g_j\leq 1$ on $\K$ for every $j$, and the family $\{g_j,1-g_j\}$ generates $\R[\x]$. If $f$ is strictly positive on $\K$ then 
\[f\,=\,\sum_{\alpha,\beta\in\N^m}c_{\alpha\beta}\,\prod_{j}g_j^{\alpha_j}\,\prod_{\ell}(1-g_\ell)^{\beta_\ell},\]
for some finitely many nonnegative scalars $(c_{\alpha\beta})$.
\end{theorem}
\subsection{Main result}

With $\K$ as in (\ref{setk}) we make the following assumption:
\begin{assumption}
\label{ass1}
$\K$ is compact and $0\leq g_j\leq 1$ on $\K$ for all
$j=1,\ldots,m$. Moreover, the family of polynomials $\{g_j, 1-g_j\}$ generates the algebra $\R[\x]$.
\end{assumption}
Notice that if $\K$ is compact and Assumption \ref{ass1} does not hold, one may always rescale
the variables $x_i$ so as to have $\K\subset [0,1]^n$, and then add redundant constraints
$0\leq x_i\leq 1$ for all $i=1,\ldots,m$. Then the family $\{g_j,1-g_j\}$ (which includes $x_j$ and $1-x_j$ for all $j$)
generates the algebra $\R[\x]$ and Assumption \ref{ass1} holds.

With $d\in\N$ and $0\leq \lambda=(\lambda_{\alpha\beta})$, $(\alpha,\beta)\in\N^{2m}_{d}$, 
let $\lambda\mapsto G_d(\lambda)$ be the function defined in (\ref{def-gedelambda}),
with associated problem:
\begin{equation}
\label{newhierarchy-def1}
\rho_d\,=\,\max_{\blambda}\:\{\:G_d(\lambda)\::\: \blambda\geq0\}.
\end{equation}
Observe that $G_d(\lambda)\leq f^*$ for all $\lambda\geq0$, and 
computing $\rho_d$ is just solving the Lagrangian relaxation of problem $\tilde{\P}_d$ in (\ref{pbtildeP}).

\begin{theorem}
\label{thmain}
Let $\K$ be as in (\ref{setk}), $f\in\R[\x]$, $d\in\N$, and let Assumption \ref{ass1} hold. Consider problem 
(\ref{newhierarchy-def1}) associated with $\P$
and with optimal value $\rho_d$. 
Then the sequence $(\rho_d)$, $d\in\N$,  is monotone nondecreasing and $\rho_d\to f^*$ as $d\to\infty$.
\end{theorem}
{\it Proof.}
We first prove that $\rho_{d+1}\geq\rho_d$ for all $d$, so that  the sequence $(\rho_d)$, $d\in\N$,  is monotone nondecreasing.
Let $0\leq\lambda=(\lambda_{\alpha\beta})$ with $(\alpha,\beta)\in\N^{2m}_{d}$. Then
$0\leq\tilde{\lambda}$ with $\tilde{\lambda}_{\alpha\beta}=\lambda_{\alpha\beta}$ whenever $(\alpha,\beta)\in\N^{2m}_{d}$, and
$\tilde{\lambda}_{\alpha\beta}=0$ whenever $\vert\alpha+\beta\vert>d$,
is such that $G_{d+1}(\tilde{\lambda})=G_{d}(\lambda)$ and so $\rho_{d+1}\geq\rho_d$.
Next, let $\epsilon>0$ be fixed, arbitrary. The polynomial $f-f^*+\epsilon$ is positive on $\K$ and therefore, by 
\cite{Stengle74}, \cite[Theorem 2.23]{lasserre-imperial},
\[f-(f^*-\epsilon)\,=\,\sum_{(\alpha,\beta)\in\N^{2m}_d}c^\epsilon_{\alpha\beta}\,\left(\prod_{j=1}^mg_j^{\alpha_j}\right)\left(\prod_{j=1}^m(1-g_j)^{\beta_j}\right),\]
for some nonnegative vector of coefficients $\c^\epsilon=(c_{\alpha\beta}^\epsilon)$. Equivalently, 
\[f-\sum_{(\alpha,\beta)\in\N^{2m}_d}c^\epsilon_{\alpha\beta}\,\left(\prod_{j=1}^mg_j^{\alpha_j}\right)\left(\prod_{j=1}^m(1-g_j)^{\beta_j}\right)\,=\,(f^*-\epsilon).\]
Letting
\[d_\epsilon\,:=\,\max_{\alpha,\beta}\:\{\vert\alpha+\beta\vert\::\:c^\epsilon_{\alpha\beta}>0\},\]
we obtain $f^*\geq G_{d_\epsilon}(\c^\epsilon)\,=\,f^*-\epsilon$. And so 
\[f^*\,\geq\,\max_\lambda\{G_{d_\epsilon}(\lambda)\,:\,\lambda\geq0\}\,\geq\,f^*-\epsilon.\]
As $\epsilon>0$ was arbitrary, the desired result follows.
$\Box$
\vspace{0.2cm}

\begin{cor}
\label{lagrangian}
Let $\K$ be as in (\ref{setk}), Assumption (\ref{ass1}) hold and let $\tilde{\P}_d$, $d\in\N$,  be as in (\ref{pbtildeP}).
Then for every $\epsilon>0$ there exists $d_\epsilon\in\N$ such that for every $d\geq d_\epsilon$,
the Lagrangian relaxation of $\tilde{\P}_d$, 
yields a lower bound $f^*-\epsilon\leq \rho_d\leq f^*$. 
\end{cor}
This follows from Theorem \ref{thmain} and the fact that computing $\rho_d$ is just solving the Lagrangian relaxation associated with
$\tilde{\P}_d$. So the interpretation of Corollary \ref{lagrangian} is that the Lagrangian relaxation
technique in non convex optimization can provide a lower bound as close as desired to the global optimum $f^*$ provided that
it is applied to an equivalent formulation of $\P$ that contains sufficiently many 
redundant constraints which are products of the original ones. It also provides a rigorous rationale for the well-known fact that
adding redundant constraints helps solve $\P$. Indeed, even though the new problems $\tilde{\P}_d$, $d\in\N$,
are all equivalent to $\P$, their Lagrangian relaxations are not equivalent to that of $\P$.

\subsection{Convex programs}

In this section, the set $\K$ is not assumed to be compact.

\begin{theorem}
\label{new-convex-case}
Let $\K$ be as in (\ref{setk})  and assume that $f$ and $-g_j$ are convex, $j=1,\ldots,m$. 
Moreover, assume that Slater's condition\footnote{Slater's condition holds for $\P$ if there exists $\x_0\in\K$ such that
$g_j(\x_0)>0$ for every $j=1,\ldots,m$.} holds and $f^*>-\infty$.

Then the hierarchy of convex relaxations (\ref{newhierarchy-def}) has finite convergence at step $d=1$, i.e., $\rho_1=f^*$,
and $\rho_1=G_1(\lambda^*)$ for some nonnegative $\lambda^*\in\R^m$.
\end{theorem}
{\it Proof.}
This is because the dual method applied to $\P$ (i.e. $\tilde{\P}_1$) converges, i.e.,
\begin{eqnarray*}
f^*&=&\max_{\lambda\geq0}\:\left\{\min_\x\:\{f(\x)-\sum_{j=1}^m\lambda_j\,g_j(\x)\}\right\}\\
&=&\max_\lambda\{G_1(\lambda)\,:\,\lambda\geq0\}\,=\,\rho_1.\end{eqnarray*}
Next, 
let $\lambda^{(n)}$ be a maximizing sequence, i.e., $G_1(\lambda^{(n)})\to f^*$ as $n\to\infty$.
Since Slater's condition holds (say at some $\x_0\in\K$), one has
\[G_1(\lambda^{(0)})\,\leq\,G_1(\lambda^{(n)})\,\leq\,f(\x_0)-\sum_{j=1}^m\lambda^{(n)}_j\,g_j(\x_0),\]
for all $n$, and so $\lambda^{(n)}_j\leq (f(\x_0)-G_1(\lambda^{(0)}))/g_j(\x_0)$ for every $j=1,\ldots,m$, and all $n\geq1$. So 
there is a subsequence $(n_k)$, $k\in\N$, 
and $\lambda^*\in\R^m_+$, such that $\lambda^{(n_k)}\to\lambda^*\geq0$ as $k\to\infty$.
Finally, let $\x\in\R^n$ be fixed, arbitrary.
From
\[G_1(\lambda^{(n_k)})\,\leq\,f(\x)-\sum_{j=1}^m\lambda^{(n_k)}_j\,g_j(\x),\qquad\forall\,k,\]
letting $k\to\infty$ yields
\[f^*\,\leq\,f(\x)-\sum_{j=1}^m\lambda^*_j\,g_j(\x).\]
As $\x\in\R^n$ was arbitrary, this proves $G_1(\lambda^*)\geq f^*$, which combined with
$G_1(\lambda^*)\leq f^*$ yields
the desired result $G_1(\lambda^*)=f^*$.
$\Box$
\vspace{0.2cm}

Observe that this does not hold for the LP-relaxations (\ref{kriv-cert}) where generically
$\theta_d<f^*$ for every $d\in\N$; see e.g. \cite{Lasserre02a,lasserre-imperial}.

\section{A parametrized hierarchy of\\ semidefinite relaxations}

Problem (\ref{newhierarchy-def1}) is convex but in general the objective function $G_d$
is non differentiable. 
Moreover, another difficulty is the computation of $G_d(\lambda)$ for each $\lambda\geq0$
since $G_d(\lambda)$ is the global optimum of the possibly non convex function
$(\x,\lambda)\mapsto L_d(\x,\lambda)$ defined in (\ref{functionL}).
So one strategy is to replace (\ref{newhierarchy-def1}) by a {\it simpler} convex problem (while preserving
the convergence property) as follows.

\subsection{Interpreting the LP-relaxations}
Observe that  the LP-relaxations (\ref{kriv-cert}) can be written
\begin{equation}
\label{lp-interpret}
\theta_d\,=\,\max_{\lambda\geq0,t}\:\left\{\:t\::\: L_d(\x,\lambda)-t\,=\,0,\quad \forall \x\in\R^n\:\right\},
\end{equation}
where $L_d$ has been defined in (\ref{functionL}).

And so the LP-relaxations (\ref{kriv-cert}) can be interpreted as simplifying (\ref{newhierarchy-def1}) by
restricting the nonnegative orthant $\{\lambda: \lambda\geq0\}$ to its subset of
$\lambda$'s that make the polynomial $\x\mapsto L(\x,\lambda)-t$ constant and equal to zero, instead of being only nonnegative. 
This subset being a polyhedron, solving (\ref{lp-interpret}) reduces to solving a linear program.
At first glance, such an a priori simple and naive brute force simplification might seem unreasonable
(to say the least). But of course the LP-relaxations (\ref{kriv-cert}) where not defined this way.
Initially, the Sherali-Adams' RLT hierarchy \cite{SheraliAdams90} was introduced for 0/1 programs and 
finite convergence was proved by using {\it ad hoc} arguments. But in fact, the rationale behind 
convergence of the more general LP-relaxations (\ref{kriv-cert}) is the Krivine-Stengle positivity certificate \cite[Theorem 2.23]{lasserre-imperial}. 

However, even though this brute force simplification
still preserves the convergence $\theta_d\to f^*$ thanks to \cite[Theorem 2.23]{lasserre-imperial}, we have already mentioned that it also implies serious theoretical (and practical) drawbacks for
the resulting LP-relaxations (like slow asymptotic convergence for convex problems and numerical ill-conditioning).

\subsection{A parametrized hierarchy of semidefinite relaxations} However, inspired by this interpretation
we propose a systematic way to improve the LP-relaxations (\ref{kriv-cert}) along the same lines
but by doing a much less brutal simplification of (\ref{newhierarchy-def1}). Indeed,
one may now impose on the same
nonnegative polynomial $\x\mapsto L(\x,\lambda)-t$ to be a
sum of squares (SOS) polynomial $\sigma$ of degree at most $2k$ (instead of being
constant and equal to zero as in (\ref{lp-interpret})), and solve the resulting hierarchy of optimization problems:

\begin{equation}
\label{lp-semi-k}
\begin{array}{rl}
q^k_d\,=\,\displaystyle\max_{\lambda,t,\sigma}&t\\
\mbox{s.t.}&L_d(\x,\lambda)-t=\sigma,\quad \forall \x\in\R^n\\
&\lambda\geq0,\quad\sigma\in\Sigma[\x]_k\end{array},
\end{equation}
with $d=1,2,\ldots$, and parametrized by $k$, fixed.
(Recall that  $\Sigma[\x]_k$ denotes the set of SOS polynomials of degree at most $2k$.) To see that (\ref{lp-semi-k}) is a semidefinite program, write
\[\x\mapsto L_d(\x,\lambda)-t\,:=\,\sum_{\beta\in\N^n_s}L_\beta(\lambda,t)\,\x^\beta,\]
where $s=d\max_j[{\rm deg}\,g_j]$ and $L_\beta(\lambda,t)$ is linear in $(\lambda,t)$ for each $\beta\in\N^n_s$.

Next, for $k\in\N$ such that $2k\leq s$, let $\v_k(\x)$ be the vector of the monomial basis $(\x^\beta)$, $\beta\in\N^n_k$, of $\R[\x]_k$, and write
\[\v_k(\x)\,\v_k(\x)^T\,=\,\sum_{\beta\in\N^n_{2k}}\x^\beta\,\B_\beta,\]
for some appropriate real symmetric matrices $(\B_\beta)$, $\beta\in\N^n_{2k}$. 
Then problem (\ref{lp-semi-k}) is the semidefinite program:
\begin{equation}
\label{sdp-k}
\begin{array}{rcll}
q^k_d\,=\,\displaystyle\max_{\lambda,t,\Q}&t&&\\
\mbox{s.t.}&L_\beta(\lambda,t)&=&\langle \B_\beta,\Q\rangle,\quad \forall \beta\in\N^n_{2k}\\
&L_\beta(\lambda,t)&=&0,\quad\forall \beta\in\N^n_s,\:\vert\beta\vert>2k\\
&\lambda&\geq&0;\:\Q=\Q^T\succeq0,
\end{array}\end{equation}
where $\Q$ is a ${n+k\choose n}\times {n+k\choose n}$ real symmetric matrix.

Of course $q^k_d\geq\theta_d\:(=q^0_d)$ for all $d$ because with $\sigma=0$ one retrieves (\ref{kriv-cert}).
Moreover in the semidefinite program (\ref{sdp-k}),
the semidefinite constraint $\Q\succeq0$ is concerned with a 
real symmetric ${n+k\choose n}\times {n+k\choose n}$ matrix, independently of the rank $d$ in the hierarchy. For instance
if $k=1$ then $\sigma$ is a quadratic SOS and $\Q$ has size $(n+1)\times (n+1)$.
In other words, even if the number of variables $\lambda=(\lambda_{\alpha\beta})$ increases fast
 with $d$, the LMI constraint $\Q\succeq0$ has fixed size, in contrast with the semidefinite relaxations
 (\ref{put-cert}) where the size of the LMIs increases with $d$. And it is a well-known fact that 
 crucial for solving semidefinite program is the size of the
 LMIs involved rather than the number of variables.

\subsection{Sherali-Adams' RLT for 0/1 programs}
\label{sherali}
Consider 0/1 programs with $f\in\R[\x]$, and feasible set $\K=\{\x:\A\x\leq\b\}\cap\{0,1\}^n$, for some
real matrix $\A\in\R^{m\times n}$ and some vector $\b\in\R^m$. The Sherali-Adams's RLT hierarchy \cite{SheraliAdams90} belongs to the family
of LP-relaxations (\ref{kriv-cert}) but with a more specific form since $\K\subset [0,1]^n$. Notice
that the family $\{1,x_1,(1-x_1),\ldots,x_n,(1-x_n)\}$ generates the algebra $\R[\x]$. Let
$g_\ell(\x)=(b-\A\x)_\ell$, $\ell=1,\ldots,m$, and $g_0(\x)=1$. 

Following the definition of the Sherali-Adams' RLT in \cite{SheraliAdams90},
the resulting linear program at step $d$ in the hierarchy reads:
\[\theta_d=\displaystyle\max_{\lambda\geq0,t,h}\left \{ t\::\:f(\x)-t\,=\,\sum_{i=1}^nh_i(\x)\,x_i(1-x_i)\right.\]
\[+\sum_{\ell=0}^m\:\displaystyle\sum_{\stackrel{I,J\subset\{1,\ldots,n\}}{
I\cap J=\emptyset; \vert I\cup J\vert \leq d}}
\lambda^\ell_{IJ}\:g_\ell(\x)\prod_{i\in I}x_i\prod_{j\in J}(1-x_j);\]
\begin{equation}
\label{sherali-def}
\left.h_i\in\R[\x]_{d-1} \quad i=1,\ldots,n\right\},\end{equation}
where $\lambda$ is the nonnegative vector $(\lambda^\ell_{IJ})$. 
(If there are linear equality constraints $g_\ell(\x)=0$ the corresponding variables $\lambda^\ell_{IJ}$ are
not required to be nonnegative.)
So all products between the $g_\ell$'s are ignored
(see the paragraph before Lemma 1 in \cite[p. 414]{SheraliAdams90}) even though they might help tighten the relaxations.
In the literature the dual LP of (\ref{sherali-def}) is described rather than (\ref{sherali-def}) itself.

In this context, the problem $\tilde{\P}_d$ equivalent to $\P$  and defined in (\ref{pbtildeP})
by adding redundant constraints formed with products of original ones, reads:
\[\begin{array}{rl}
\min\{f(\x)\::& \x^\alpha \,x_j(1-x_j)\,=\,0,\: j=1,\ldots,n;\:\alpha\in\N^n_{d-1};\\
&g_\ell(\x)\displaystyle\prod_{i\in I}x_j\displaystyle\prod_{j\in J}(1-x_j)\geq0,\quad \ell=0,\ldots,m,\\
&I,J\subset\{1,\ldots,n\};\:I\cap J=\emptyset;\:\vert I\cup J\vert\leq d\}.
\end{array}\]

Hence the 0/1 analogue of (\ref{lp-semi-k})  reads

\[q^k_d=\displaystyle\max_{\lambda\geq0,t,h}\left \{ t\::\:f(\x)-t\,=\,\sigma(\x)+\sum_{i=1}^nh_i(\x)\,x_i(1-x_i)\right.\]
\[+\sum_{\ell=0}^m\:\displaystyle\sum_{\stackrel{I,J\subset\{1,\ldots,n\}}{
I\cap J=\emptyset; \vert I\cup J\vert \leq d}}
\lambda^\ell_{IJ}\:g_\ell(\x)\prod_{i\in I}x_i\prod_{j\in J}(1-x_j);\]
\begin{equation}
 \label{lp-semi-k-01}
\left.\sigma\in\Sigma[\x]_{k};\quad h_i\in\R[\x]_{d-1} \quad i=1,\ldots,n\right\}.\end{equation}

For 0/1 programs with linear or quadratic objective function, and for every $k\geq1$,
the first semidefinite relaxation (\ref{lp-semi-k-01}), i.e., with $d=2$,
is at least as powerful as that of the standard hierarchy of semidefinite relaxations (\ref{put-cert}).
Indeed (\ref{lp-semi-k-01}) contains products $g_\ell(\x)x_j$ or $g_\ell(\x)(1-x_k)$, for all $(\ell,j,k)$, which do  to not appear in (\ref{put-cert}) with $d=1$.
And so in particular, the first such relaxation for MAXCUT has the celebrated Goemans-Williamson's performance guarantee
while the standard LP-relaxations (\ref{kriv-cert}) do not. On the other hand, for 0/1 problems
and for the parameter value $k=1$,
the hierarchy (\ref{lp-semi-k-01}) is what is called the {\it Sherali-Adams + SDP} hierarchy 
(basic SDP-relaxation + RLT hierarchy)
in e.g. Benabas and Magen \cite{benabbas1} and Benabbas et al. \cite{benabbas2}; 
and in \cite{benabbas1,benabbas2} the authors show that any (constant) level $d$ of this hierarchy,
viewed as a strengthening of the basic SDP-relaxation, does not make the integrality gap decrease.

In fact, and in view of our previous analysis, the ``Sherali-Adams + SDP" hierarchy should be viewed as 
a {\it (level $k=1$)}-strengthening of the basic Sherali-Adams'  LP-hierarchy
(\ref{sherali-def})  rather than a strengthening of the basic SDP relaxation.

 \section{Comparing with standard \\LP-relaxations}
 
 As asked in introduction:\\ {\it What do we gain by going from the LP hierarchy (\ref{kriv-cert}) to the 
 semidefinite hierarchy (\ref{sdp-k}) parametrized by $k$?}
 Some answers are provided below.

 \subsection{Convex problems}
 
 Recall that a highly desirable property for a general purpose method 
aiming at solving NP-hard optimization problems, is to behave efficiently when 
applied to a class of problems considered relatively easy to solve. Otherwise one might raise reasonable doubts on  its
efficiency for more difficult problems not only in a worst-case sense but also in average. And convex problems $\P$ as in (\ref{def-pb})-(\ref{setk}), i.e., when $f,-g_j$ are convex,
form the most natural class of problems which are considered easy to solve by some standard methods of Non Linear Programming.

 \begin{theorem}
 \label{th-final}
With $\P$ as in (\ref{def-pb})-(\ref{setk}) let $f,-g_j$ be convex, $j=1,\ldots,m$, let Slater's condition hold and let $f^*>-\infty$.
Then:

(a) If $\max[{\rm deg}\,f,\:{\rm deg}\,g_j]\leq 2$ then $q^1_1=f^*$, i.e., the first relaxation of the hierarchy 
(\ref{lp-semi-k}) parametrized by $k=1$, is exact.

(a) If $\max[{\rm deg}\,f,\:{\rm deg}\,g_j]\leq 2k$ and $f,-g_j$ are all SOS-convex, then $q^k_1=f^*$, i.e., the first relaxation of the hierarchy 
(\ref{lp-semi-k}) parametrized by $k$, is exact.
 \end{theorem}
 {\it Proof.}
 Under the assumptions of Theorem \ref{th-final}, $\P$ has 
 a minimizer $\x^*\in\K$ and the Karush-Kuhn-Tucker optimality conditions hold
 at  $(\x^*,\lambda^*)\in\K\times\R^m_+$ for some $\lambda^*\in\R^m_+$.
 And so if $k=1$, the Lagrangian polynomial  $L_1(\cdot,\lambda^*)-f^*$ is a nonnegative quadratic polynomial
 and so an SOS $\sigma^*\in\Sigma[\x]_1$. Therefore
 as $q^1_d\leq f^*$ for all $d$,  the triplet $(\lambda^*,f^*,\sigma^*)$ is an optimal solution of (\ref{lp-semi-k}) with $k=d=1$, which proves (a).
 
 Next, if $k>1$ and $f,-g_j$ are all SOS-convex then so is the Lagrangian polynomial
 $L_1(\cdot,\lambda^*)-f^*$. In addition, as 
 $\nabla_\x L_1(\x^*,\lambda^*)=0$ and $L_1(\x^*,\lambda^*)-f^*=0$,
 the polynomial $L_1(\cdot,\lambda^*)-f^*$ is SOS; see e.g. Helton and Nie \cite[Lemma 4.2]{nie-helton-survey}. Hence 
 $L_1(\cdot,\lambda^*)-f^*=\sigma^*$ for some $\sigma^*\in\Sigma[\x]_k$, and again,
 the triplet $(\lambda^*,f^*,\sigma^*)$ is an optimal solution of (\ref{lp-semi-k}) with $d=1$, which proves (b).
 $\Box$
 \vspace{0.2cm}
 
 Hence  by simplifying (\ref{newhierarchy-def}) in a less brutal manner than in (\ref{kriv-cert})
 one recovers a nice and highly desirable property for the resulting hierarchy. The price to pay is to pass from solving a hierarchy of LPs to solving  hierarchy of semidefinite programs; however the increase in complexity is controlled
 by the parameter $k$ since the size of the LMI in the semidefinite program (\ref{sdp-k}) is $O(n^k)$, independently of the rank $d$ in the hierarchy.
 
\subsection{Obstructions to Exactness}
On the other hand, for non convex problems, exactness at level-$d$ of the hierarchy (\ref{lp-semi-k}), i.e., finite convergence
after $d$ rounds, still implies restrictive conditions on the problem:

\begin{cor}
\label{badnews}
Let $\P$ be as in (\ref{def-pb})-(\ref{setk}) and let Assumption \ref{ass1} hold.
Let $\x^*\in\K$ be a global minimizer and let 
$I_1(\x^*):=\{j\in\{1,\ldots,m\}:\,g_j(\x^*)=0\}$ and
$I_2(\x^*):=\{j\in\{1,\ldots,m\}:\,(1-g_j(\x^*))=0\}$ 
be the set of active constraints at $\x^*$. Let $0\leq k\in\N$ be fixed.

The level-$d$ semidefinite relaxation (\ref{lp-semi-k}) is exact
only if $f^*$ (resp. $\x^*\in\K$) is also the global optimum (resp. a global minimizer) for the problem
\begin{equation}
\label{pb-variety}
\min_\x\:\{f(\x)\::\: \x\in \V\},\end{equation}
where $\V\subset\R^n$ (see (\ref{real-variety}) below) is a variety defined from 
some products of the polynomials $g_j$'s and $(1-g_j)$'s. 
And if $k=0$ then $f$ must be constant on the variety $\V$.
\end{cor}
{\it Proof.}
If (\ref{lp-semi-k}) is exact at level $d\in\N$, then
\[L_d(\x,\lambda^*)-f^*=\sigma(\x),\qquad\forall \x\in\R^n,\]
for some $\lambda^*\geq0$ and some $\sigma\in\Sigma[\x]_k$. Equivalently,
\begin{eqnarray*}
f(\x)-f^*\,=\,\sigma(\x)&+&\sum_{(\alpha,\beta)\in\N^{2m}_{d}}\lambda^*_{\alpha\beta}\left(\prod_{j=1}^mg_j(\x)^{\alpha_j}\right)\times\\
&&\left(\prod_{j=1}^m(1-g_j(\x))^{\beta_j}\right).\end{eqnarray*}
Then evaluating at $\x=\x^*$ yields $\sigma(\x^*)=0$ and 
\[\lambda^*_{\alpha\beta}\,>0\,\Rightarrow\quad\left\{\begin{array}{l}
\exists j\in I_1(\x^*)\mbox{ s.t. }\alpha_j>0,\:\mbox{or}\\
\exists j\in I_2(\x^*)\mbox{ s.t. }\beta_j>0.\end{array}\right.\]
So let $\Omega\,:=\,\{(\alpha,\beta)\in\N^{2m}_{d}\::\:\lambda^*_{\alpha\beta}>0\}$ and for every
$(\alpha,\beta)\in\Omega$ let
\begin{eqnarray*}
J^1_{\alpha\beta}&:=&\,\{j\in I_1(\x^*)\,:\,\alpha_j>0\},\\
J^2_{\alpha\beta}&:=&\{j\in I_2(\x^*)\,:\,\beta_j>0\}.\end{eqnarray*}
Next, define $\V\subset\R^n$ to be the real variety:
\begin{eqnarray}
\nonumber
\left\{\x\in\R^n\right.&:&
\left(\prod_{j\in J^1_{\alpha\beta}}g_j(\x)\right)
\left(\prod_{j\in J^2_{\alpha\beta}}(1-g_j(\x))\right)\,=\,0,\\
\label{real-variety}
&&\left.\forall\:(\alpha,\beta)\in\Omega\:\right\}.\end{eqnarray}
Then for every $\x\in\V$, one obtains $f(\x)-f^*=\sigma(\x)\geq0$, which means that
$f^*$ is the global minimum of $f$ on $\V$. If $k=0$ then $\sigma$ is constant and equal to zero.
And so $f(\x)-f^*=0$ for all $\x\in\V$.
$\Box$
\vspace{0.2cm}

Hence Corollary \ref{badnews} shows that
exactness at some step $d$ of the hierarchy (\ref{lp-semi-k}) imposes rather restrictive conditions on problem $\P$.
Namely, the global optimum $f^*$ (resp. the global minimizer $\x^*\in\K$) must also be the global optimum
(resp. a global minimizer) of problem (\ref{pb-variety}). For instance, suppose that
only one constraint, say $g_k(\x)\geq0$, is active at $\x^*$. Then $f^*$ (resp. $\x^*$) is also the global minimum 
(resp. a global minimizer) of $f$ on the variety 
$\{\x\,:\,g_k(\x)=0\}$.  And if $k=0$ then $f$ must be constant on the variety $\V$!

\begin{example} If $\K$ is the (compact) polytope $\{\x: \a_j^T \x\leq 1,\,j=1,\ldots,m\}$ for some vectors
$(\a_j)\subset\R^n$, then invoking a result by Handelman \cite{Handelman88}, one does not need the polynomials $\{1-g_j\}$ in the definition
(\ref{functionL}) of $L_d$. So for instance, suppose that $I_1(\x^*)=\{\ell\}$
at a global minimizer $\x^*\in\K$.
Then exactness at some step $d$ of the hierarchy (\ref{lp-semi-k}) imposes that $f^*$ should also be the global minimum of
$f$ on the whole hyperplane $\V=\{\x: \a_\ell^T\x=1\}$; for non convex functions $f$,
this is a serious restriction. Moreover, if $k=0$ then $f$ must be constant on the hyperplane $\V$. 
\end{example}

Concerning exactness for 0/1 polynomial optimization:
\begin{cor}
\label{sherali-exact}
Let $\K=\{\x:\A\x\leq\b\}\cap\{0,1\}^n$ and let
$\x^*\in\K$ be an optimal solution of $f^*=\min \{f(\x):\x\in\K\}$.
Assume that $\A\x^*<\b$, i.e., no constraint is active at $\x^*$. 

(a) The Sherali-Adams' RLT relaxation (\ref{sherali-def}) is exact at step
$d$ in the hierarchy only if $f(\x)=f(\x^*)=f^*$ for all
$\x$ in the set
\[\V:=\{\x\in\{0,1\}^n\::\:\prod_{i\in I}x_i\prod_{j\in J}(1-x_j)=0,\quad (I,J)\in Q\},\]
where $Q$ is some finite set of couples $(I,J)$ satisfying
$I\cap J=\emptyset$ and $\vert I\cup J\vert\leq d$. 

(b)  Similarly, the semidefinite relaxation (\ref{lp-semi-k-01}) is exact at step $d$ only if 
$\x^*$ is also a global minimizer of $\min \{f(\x)\,:\:\x\in\V\}$ for some $\V$
as in (a).
\end{cor}
{\it Proof.}
(a) Exactness implies that the polynomial $\x\mapsto f(\x)-f^*$ has 
the representation described in (\ref{sherali-def}) for some polynomials $(h_i)\subset\R[\x]_{d-2}$ and some nonnegative scalars 
$(\lambda^\ell_{IJ})$.
Evaluating both sides of (\ref{sherali-def}) at $\x=\x^*$ and using $g_\ell(\x^*)>0$ for all $\ell=0,\ldots,m$, yields
\begin{equation}
\label{aux-sherali}
\lambda_{IJ}^\ell>0\quad\Longrightarrow\quad\prod_{i\in I}x^*_i\prod_{j\in J}(1-x^*_j)=0.\end{equation}
Let $\V$ be as in Corollary \ref{sherali-exact} with $Q:=\{(I,J)\::\:\mbox{$\exists\,\ell$ s.t. }\lambda^\ell_{IJ}>0\}$. Then
from the representation of $\x\mapsto f(\x)-f^*$ in (\ref{sherali-def}) we obtain
$\f(\x)-f^*=0$ for all $\x\in\V$ and the result follows. For (b) a similar argument is valid but now using the representation
of $f(\x)-f^*$ described in (\ref{lp-semi-k-01}). And so exactness yields (\ref{aux-sherali}) as well as
 $\sigma(\x^*)=0$. Next, for every $\x\in \V$ we now obtain $f(\x)-f^*=\sigma(\x)\geq0$ because $\sigma$ is SOS.
$\Box$
\vspace{0.2cm}

The constraints $\A\x\leq\b$ play no explicit role in the definition
of the set $\V$. Moreover, if $f$ discriminates all points of the hypercube $\{0,1\}^n$
then exactness of the Sherali-Adams' RLT
implies that $\V$ must be the singleton $\{\x^*\}$.

\subsection*{On the hierarchy of semidefinite relaxations} 

Similarly, the hierarchy of semidefinite relaxations (\ref{put-cert}) also has an interpretation 
in terms of simplifying an {\it extended} Lagrangian relaxation of $\P$. 
Indeed consider the hierarchy of optimization problems
\begin{eqnarray}
\nonumber
\omega_d\,:=\,\max_{\sigma_j}\:\{H(\sigma_1,\ldots,\sigma_m)&:&
{\rm deg}(\sigma_j\,g_j)\leq 2d,\sigma_j\in\Sigma[\x]\\
\label{ex-lagrangian}
&&j=1,\ldots,m\},
\end{eqnarray}
$d\in\N$, where $\sigma\mapsto H(\sigma_1,\ldots,\sigma_m)$ is the function
\[H(\sigma_1,\ldots,\sigma_m)\,:=\,\min_\x\:\{f(\x)-\sum_{j=1}^m\sigma_j(\x)\,g_j(\x)\:\}.\]
For each $d\in\N$, problem (\ref{ex-lagrangian}) is an obvious relaxation of $\P$ and in fact 
is an {\it extended Lagrangian relaxation} of $\P$ where the multipliers are now allowed to be SOS 
polynomials with a degree bound, instead of constant nonnegative polynomials (i.e., SOS polynomials of degree zero).

If $\K$ is compact and the quadratic module 
\[Q(g)\,:=\,\{\sum_{j=0}^m\sigma_j\,g_j\::\:\sigma_j\in\Sigma[\x],\quad j=0,1,\ldots,m\}\]
(where $g_0=1$) is Archimedean, then $\omega_d\to f^*$ as $d\to\infty$. But of course,
and like for the usual Lagrangian, minimizing the extended Lagrangian
\[\x\mapsto L(\x,\sigma):=f(\x)-\sum_{j=1}^m\sigma_j(\x)\,g_j(\x),\]
 is in general an NP-hard problem. In fact, writing (\ref{ex-lagrangian}) as
 \begin{eqnarray*}
 \omega_d=\max_{t,\sigma}\{t&:&f(\x)-\sum_{j=1}^m\sigma_j(\x)\,g_j(\x)-t\,\geq\,0\quad\forall\x\,;\\
 &&{\rm deg}(\sigma_jg_j)\leq 2d\},\end{eqnarray*}
 the semidefinite relaxations (\ref{put-cert}) simplify
 (\ref{ex-lagrangian}) by imposing on the nonnegative polynomial
 $\x\mapsto f(\x)-\sum_j\sigma_j(\x)g_j(\x)-t$ to be an SOS polynomial $\sigma_0\in\Sigma[\x]_d$
 (rather than just being nonnegative).
 
 But the spirit is different from the LP-relaxations as there is no problem $\tilde{\P}_d$ obtained from $\P$
 by adding
finitely many  redundant constraints and equivalent to $\P$. Instead of adding more and more redundant constraints
 and doing a standard Lagrangian relaxation to $\tilde{\P}_d$, one applies an extended Lagrangian relaxation to 
 $\P$ with SOS multipliers of increasing degree (instead of nonnegative scalars). And in contrast to
 LP-relaxations, there is no obstruction to exactness (i.e., finite convergence). In fact, it is quite the opposite since
 as demonstrated recently in Nie \cite{nie2012}, finite convergence is generic!

 \section{Conclusion}
  
We have shown that the hierarchy of LP-relaxations (\ref{kriv-cert}) has a rather surprising 
interpretation in terms of the Lagrangian relaxation applied to
a problem $\tilde{\P}$ equivalent to $\P$ (but with redundant
constraints formed with product of polynomials defining the original constraints of $\P$).
Indeed it consists of the brute force simplification of imposing 
on a certain nonnegative polynomial to be the constant polynomial equal to zero, a very restrictive
condition.

However, inspired by this interpretation,  one has provided a systematic strategy
to improve the LP-hierarchy
by doing a much less brutal simplification. That is,
one now imposes on the same nonnegative polynomial to be an SOS polynomial 
whose degree $k$ is fixed in advance and parametrizes the whole hierarchy. Each convex relaxation is now a semidefinite program but 
whose LMI constraint has fixed size $O(n^k)$. Hence, the resulting families of parametrized relaxations achieve a compromise
between the hierarchy of semidefinite relaxations (\ref{put-cert}) limited to 
problems of modest size and the LP-relaxations (\ref{kriv-cert}) that theoretically can handle problems of 
larger size but with a poor behavior when applied to convex problems.

\end{document}